\documentclass[journal]{IEEEtran}

\ifCLASSINFOpdf
\else
\fi
\usepackage{graphics} 
\usepackage{epsfig} 
\usepackage{amsmath} 
\usepackage{amssymb}  
\usepackage{amsthm}

\usepackage{bm}
\usepackage{multirow}
\usepackage{tabularx}
\usepackage{algpseudocode}
\usepackage{algorithm}
\usepackage{epstopdf}
\usepackage{circuitikz}
\usepackage{cite, balance}

\allowdisplaybreaks

\DeclareMathOperator{\diag}{diag}

\theoremstyle{definition}

\begin{document}
%
\title{Counterexample to Equivalent Nodal Analysis for Voltage Stability Assessment}
%
%
%

\author{Bai~Cui,~\IEEEmembership{Member,~IEEE} and Zhaoyu Wang, \IEEEmembership{Member, IEEE}
\thanks{B. Cui is with the Energy Systems Division, Argonne National Laboratory, 9700 Cass Avenue, Lemont, IL 60439 (e-mail: bcui@anl.gov).}
\thanks{Z. Wang is with the Department of Electrical and Computer Engineering, Iowa State University, Ames, IA 50011 USA (e-mail: wzy@iastate.edu).}
}

\maketitle

\begin{abstract}
Existing literature claims that the $L$-index for voltage instability detection is inaccurate and proposes an improved index quantifying voltage stability through system equivalencing. The proposed stability condition is claimed to be exact in determining voltage instability. We show the condition is incorrect through simple arguments accompanied by demonstration on a two-bus system counterexample. 
\end{abstract}


%
\IEEEpeerreviewmaketitle

\section{Introduction} \label{sect1}

\IEEEPARstart{P}{ower} flow equations are ubiquitous in power system analysis. It is shown in \cite{Sauer90} that the singularity of power flow Jacobian is closely related to the loss of steady-state voltage stability. In power system on-line monitoring and control, the knowledge of stability margins of the current operating point, as well as that under system parameter variations and plausible contingencies, are vital for maintaining a steady and reliable operation. Simple on-line indicators/indices to characterize the relative distance to voltage stability boundary as well as to determine the onset of voltage instability have been proposed for both local and wide-area monitoring, with little computational overhead \cite{Simpson-Porco16, Bolognani16, WangC16, Kessel86, WangY11, Xu12, Pordanjani13, Wang13, Vu99, Liu14, Cui17}.

\begin{figure}[!t]
	\centering
	\begin{circuitikz}
		\draw[thick] (0,0) node[oscillator]{};
		\draw[very thick] (0,0) -- (1,0);
		\draw[very thick] (1,-0.5) -- (1, 0.5);
		\node [above] at (1,0.5) {$E_{eq}$};
		\draw[thick] (1,0) to[generic] (4,0);
		\draw[very thick] (1,0) -- (1.92,0);
		\draw[very thick] (3.04,0) -- (4,0);
		\node [above] at (2.5, 0.3) {$Z_{eq}$};
		\draw[very thick] (4,-0.5) -- (4, 0.5);
		\node [above] at (4,0.5) {$V_L$};
		\draw[very thick, ->] (4,0) -- (5,0);
		\node [right] at (5,0) {$S_L$};
		
	\end{circuitikz}
	\caption{Equivalent two-bus system model with equivalent generator voltage $E_{eq}$ and equivalent transmission line with impedance $Z_{eq}$.}
	\label{fig1}
\end{figure}
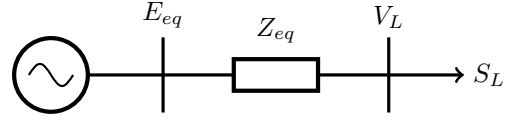

These indices can be categorized into parameter-based \cite{Simpson-Porco16, Bolognani16, WangC16} and state-based \cite{Kessel86, WangY11, Xu12, Pordanjani13, Wang13, Vu99, Liu14, Cui17} ones. Parameter-based indices characterize inner approximations of the power flow solvability boundary in power injection space beyond which no high-voltage power flow solutions exist. These indices provide rigorous lower bounds of maximum allowable injections such that voltage stability can be maintained. Active research efforts are currently undertaken to improve the conservativeness issue \cite{Dorfler18}. State-based indices quantify the system stability level through bus voltage measurements, and can be further categorized into model-based \cite{Kessel86, WangY11, Xu12, Pordanjani13, Wang13} and measurement-based \cite{Vu99, Liu14, Cui17} methods. These model-based indices are based on the idea of generalizing the classic `impedance matching' condition for a two-bus system \cite{Vournas15}, where the entire system seen from the load bus is modeled as an equivalent voltage source connected to the local load through an equivalent transmission line, as shown in Fig. \ref{fig1}. We note that special care should be taken to extend the impedance matching condition to cases where the equivalent voltage and impedance are not constant. Under certain assumptions on system modeling and loading patterns, the heuristic two-bus approximation works remarkably well in quantifying the system stability level, as demonstrated in \cite{Kessel86, WangY11}. However, it should be noted that because of the heuristic nature of these models, impedance matching between the equivalent line impedance and the local load impedance does not necessarily coincide with the point of voltage instability \cite{Cui17}. This is to be contrasted with measurement-based indices proposed in \cite{Vu99, Liu14, Cui17}, where local voltage to power sensitivities are captured and used to construct an equivalent two-bus system. In these models, the impedance matching occurs at the point of voltage instability as a result of infinite voltage to power sensitivities at the point of voltage instability. 

As a state- and model-based index, the method proposed in \cite{Lee16} is claimed to be able to determine the exact steady-state voltage instability point where load power injections cannot be further increased. We show that this claim is incorrect through simple arguments and comparison with a classic model-based index, as well as demonstration on a two-bus system. 

The remainder of the paper is organized as follows: we briefly introduce $L$-index, the classic model-based index, and our recent interpretation of the index as a power flow Jacobian singularity indicator \cite{Wang17} in Section \ref{sect:lindex}. We then introduce the acclaimed exact condition to detect voltage instability proposed in \cite{Lee16} in Section \ref{sect:ENA} and point out its problem. A two-bus system counterexample is presented in Section \ref{sect:example} to demonstrate the incorrectness of the index. Section \ref{sect:conclusion} concludes the paper.

\section{$L$-Index as a Jacobian Singularity Indicator} \label{sect:lindex}

Given bus admittance matrix $Y$ and partition it based on generator and load buses, Ohm's law dictates the following relationship between bus current and voltage vectors:
\begin{equation} \label{admit}
\begin{bmatrix}
I_G \\
I_L
\end{bmatrix} = 
\begin{bmatrix}
Y_{GG} & Y_{GL} \\
Y_{LG} & Y_{LL}
\end{bmatrix}
\begin{bmatrix}
V_G \\
V_L
\end{bmatrix},
\end{equation}
where subscripts $G$ and $L$ denote generator and load buses, respectively. Solving for $V_L$ in \eqref{admit} yields
\begin{equation}
V_L = -Y_{LL}^{-1}Y_{LG}V_G + Y_{LL}^{-1}I_L.
\label{equi1}
\end{equation}
Denote the vector of equivalent voltage as $E := -Y_{LL}^{-1} Y_{LG} V_G$ and the impedance matrix as $Z := Y_{LL}^{-1}$ (the invertibility of $Y_{LL}$ for general power system has been shown in \cite{WangC16}), \eqref{equi1} can be rewritten as 
\begin{equation} \label{eq:nodekvl}
V_L = E + ZI_L.
\end{equation}
For each load bus $i$, \eqref{eq:nodekvl} can be interpreted as a two-bus equivalent seen from the local load bus, where $E_i$ is the equivalent generator voltage and $(z_i^\top I_L)/I_{L,i}$ is the equivalent line impedance ($z_i^\top$ denotes the $i$th row of $Z$). The maximum loadability at the receiving end of a two-bus system can be solved in closed-form: both the load power and voltage of the receiving bus at the verge of voltage instability can be expressed analytically as functions of system parameters. Therefore, under the assumptions that for all load buses $i$, both $E_i$ and $(z_i^\top I_L)/I_{L,i}$ are constant while load current $I_{L,i}$ is varied, the $L$-index
\begin{equation}
	L := \max_i \left| \frac{E_i - V_i}{V_i} \right|
\end{equation}
becomes one when voltage instability occurs. Note that this is exactly the `impedance matching' condition for two-bus system \cite{Vournas15}.

The effectiveness of $L$-index has been rigorously justified in our recent paper \cite{Wang17}, which shows that $L$-index is related to the diagonal dominance of a similarity transformed power flow Jacobian matrix. Specifically, we show that when generator voltages are constant, the following matrix can be obtained through similarity transformation from a complex power flow Jacobian matrix, so that the eigenvalues (and their algebraic multiplicities) of the complex power flow Jacobian are preserved:
\begin{equation}
	J = \diag([V_i; V_i^*]) + \begin{bmatrix} 0 & B \\ B^* & 0 \end{bmatrix},
\end{equation}
where $[B_{ij}] = Z_{ij}I_j$. It follows that power flow Jacobian is nonsingular as long as $J$ is strictly diagonally dominant, which holds as long as 
\begin{equation}
	|V_i| > \| \diag(z_i^\top) I_L \|_1
\end{equation}
for all load bus $i$. By the relationship $E_i - V_i = z_i^\top I_L$, it's clear that $L$-index is closely related to the necessary condition for power flow Jacobian non-singularity. Notice, however, that based on the above discussion, a unity $L$-index is neither necessary nor sufficient for voltage instability, and deriving an exact algebraic characterization of the singularity of power flow Jabobian seems quite challenging.

\section{Equivalent Nodal Analysis} \label{sect:ENA}

Contrary to $L$-index or similar heuristic model-based voltage stability indices, the method proposed in \cite{Lee16} is claimed to be able to exactly determine the onset of voltage instability based on system parameters and nodal measurements. The method is very similar to $L$-index in spirit, both casting the system as multi-port two-bus equivalents, and then deriving stability conditions based on the two-bus impedance matching condition. In this section, we briefly review the method in \cite{Lee16}, and then point out its flaw.

Multiply both sides of \eqref{admit} by $Y^{-1}$, we obtain
\begin{equation}
V = \tilde ZI,
\end{equation}
where the impedance matrix $\tilde Z := Y^{-1}$ is the inverse of bus admittance matrix. There is already an issue here, since $Y$ is not invertible unless the system has shunt elements \cite{Kettner18}. We proceed assuming $Y$ is indeed invertible. Then, for a specific load bus $i$, we have
\begin{equation}
V_i = \sum_{j \in \mathcal{N}_G} \tilde Z_{ij}I_j + \tilde Z_{ii}\sum_{j \in \mathcal{N}_L} \frac{\tilde Z_{ij}}{\tilde Z_{ii}}I_j,
\end{equation}
where $\mathcal{N}_G$ and $\mathcal{N}_L$ stands for the index set of generator and load buses, respectively. The equivalent generator voltage, equivalent line impedance and equivalent load seen by bus $i$ are then defined as
\begin{equation}
V_{si} := \sum_{j \in \mathcal{N}_G} \tilde Z_{ij}I_j, Z_{eq,i} := \tilde Z_{ii}, S_{eq,i} := V_i\left( \sum_{j \in \mathcal{N}_L} \frac{\tilde Z_{ij}}{\tilde Z_{ii}}I_j \right)^*,
\end{equation}
respectively. For brevity, we drop subscript $i$ in the sequel with the understanding that we focus on the two-bus equivalent system seen from bus $i$. 

Let the real and imaginary parts of $\tilde Z_{eq}$ be $R_{eq}$ and $X_{eq}$. The following relationship can be derived by multiplying both sides of $S_{eq}^* = V^* (V_{s} - V)/Z_{eq} := P_{eq} + \mathrm{i}Q_{eq}$ by their complex conjugates and rearrange terms:
\begin{multline} \label{eq:quadratic}
|V|^4 + \left( 2(P_{eq}R_{eq} + Q_{eq}X_{eq}) - |V_s|^2 \right) |V|^2\\
 + |S_{eq}|^2|Z_{eq}|^2 = 0.
\end{multline}
The paper \cite{Lee16} proposes the following way to exactly identify the stability boundary: By treating \eqref{eq:quadratic} as a quadratic equation in $|V|^2$, the discriminant is given by
\begin{equation} \label{eq:discriminant}
\Delta := (|V_s|^2 - \alpha_1)(|V_s|^2 + \alpha_2),
\end{equation}
where $\alpha_1 := 2|Z_{eq}||S_{eq}|(1+\cos\phi)$, $\alpha_2 := 2|Z_{eq}||S_{eq}|$ $ (1 - \cos\phi)$, and $\phi$ is the angle difference between load power $S_{eq}$ and equivalent impedance $Z_{eq}$.

The paper \cite{Lee16} claims that \eqref{eq:quadratic} has real solutions only when $|V_s|^2 - \alpha_1 \geq 0$, and the inequality implies $|V_s|^2 = \alpha_1$ on voltage stability boundary. By setting $\Delta = 0$ and solving for $|V|^2$, the voltage magnitude at bus $i$ at voltage stability boundary is $|V|^2 = |Z_{eq}| |S_{eq}^{\max}|$, and the load impedance at stability boundary is
\begin{equation}
|Z_L| := \frac{|V|^2}{|S_{eq}^{\max}|} = |Z_{eq}|.
\end{equation}

A fundamental flaw of the proposed approach in \cite{Lee16} lies in the fact that \eqref{eq:quadratic} cannot be treated as a quadratic equation of $|V|^2$. While $P_{eq}$, $Q_{eq}$, $R_{eq}$ and $X_{eq}$ are given system parameters, $V_s$ is an implicit function of $V$, so the above analysis based on the implicit assumption that the coefficients of the quadratic equation are constant is problematic. We note that the derivations of $L$-index and its more recent extensions \cite{WangY11, Xu12, Pordanjani13, Wang13} employ similar arguments, but construct the two-bus equivalent system in such a way that the equivalent system parameters ($E_i$ and $(z_i^\top I_L)/I_{L,i}$) are approximately constant while load powers are varied, and the approximation nature of their approaches is explicitly acknowledged. On the contrary, the equivalent generator voltage $V_s$ for the proposed model in \cite{Lee16} is a weighted sum of generator currents, which is far from constant as loads ramp up.

%

The aforementioned flaw appears in Equation (18) in \cite{Lee16}, which invalidates all subsequent discussions.

\section{Numerical Counterexample} \label{sect:example}

In this section, we present a two-bus system counterexample to the proposed index in \cite{Lee16}. We show that neither $\Delta = 0$ nor $|Z_L| = |Z_{eq}|$ at the point of voltage instability.


\begin{figure}[!t]
	\centering
	\begin{circuitikz}
		\draw[thick] (0,0) node[oscillator]{};
		\draw[very thick] (0,0) -- (1,0);
		\draw[very thick] (1,-0.5) -- (1, 0.5);
		\node [above] at (1,0.5) {\small $V_1 = 1$};
		\draw[thick] (1,0) to[generic] (4,0);
		\draw[very thick] (1,0) -- (1.92,0);
		\draw[very thick] (3.04,0) -- (4,0);
		\node [above] at (2.5, 0.3) {\small $Y_{12} = -3\mathrm{i}$};
		\draw[very thick] (4,-0.5) -- (4, 0.5);
		\node [above] at (4,0.5) {\small $V_2$};
		\draw[very thick, ->] (4,0) -- (5,0);
		\node [right] at (5,0) {\small $S = P + \mathrm{i}Q$};
		
		\draw[very thick] (1,-0.25) -- (1.4, -0.25);
		\draw[thick] (1.4,-0.25) to[/tikz/circuitikz/bipoles/length=0.8cm, capacitor] (1.4,-0.9);
		\draw[thick] (1.4,-0.7) node[/tikz/circuitikz/bipoles/length=1cm, ground]{};
		\node [right] at (1.6,-0.6) {\small $Y_{sh} = \mathrm{i}$};
		
	\end{circuitikz}
	\caption{Example two-bus system model with generator voltage $V_1$, transmission line impedance $Z = 1/Y_{12}$, shunt admittance $Y_{sh}$, load bus voltage $V_2$, and load power $S = P + \mathrm{i}Q$.}
	\label{fig2}
\end{figure}
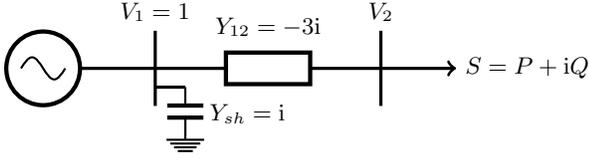

We consider a two-bus system as shown in Fig. \ref{fig2} where all quantities are specified in p.u.. Bus 1 is a slack bus with fixed voltage $V_1 = 1$ and bus 2 is a constant-power load with base power $S = 0.2\mathrm{i}$. The admittance of the line between buses 1 and 2 is $Y_{12} = -3\mathrm{i}$. To ensure the invertibility of the admittance matrix, we assume bus 1 has shunt admittance $Y_{sh} = \mathrm{i}$. The admittance matrix is
\begin{equation*}
Y = \begin{bmatrix}
Y_{12} + Y_{sh} & -Y_{12} \\ -Y_{12} & Y_{12}
\end{bmatrix}.
\end{equation*}
The inverse of the admittance matrix is
\begin{equation*}
\tilde Z := Y^{-1} = \begin{bmatrix}
1/Y_{sh} & 1/Y_{sh} \\
1/Y_{sh} & 1/Y_{sh} + 1/Y_{12}
\end{bmatrix}.
\end{equation*}

The paper \cite{Lee16} claims that the system is at its voltage stability limit when the following condition holds
\begin{equation*}
|Z_{L}| = |\tilde Z_{22}| = |1/Y_{sh} + 1/Y_{12}|,
\end{equation*} 
where $Z_L$ is the load impedance. However, it is well-known that the impedance matching condition
\begin{equation*}
|Z_{L}| = |1/Y_{12}|
\end{equation*}
holds at the stability limit for the two-bus system. Therefore, the acclaimed condition is clearly incorrect for the two-bus system.

We present the variation of $\Delta$ associated with bus 2 as load power $P+\mathrm{i}Q$ is scaled up until power flow fails to converge. It is clear from the figure that the condition $\Delta = 0$ does not necessarily correspond to the point of voltage instability.

\begin{figure}[!t]
	\centering
	\includegraphics[width=2.7in]{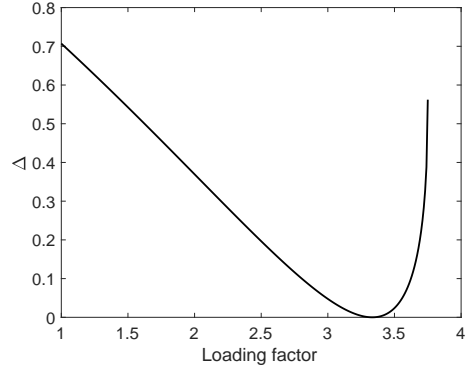}
	\caption{Variation of $\Delta$ for two-bus system as load at bus 2 is scaled up until voltage instability. System parameters are $Y_{sh} = \mathrm{i}$, $Y_{12} = -3\mathrm{i}$.}
	\label{fig:sim}
\end{figure}

\section{Conclusion} \label{sect:conclusion}

We have shown that the acclaimed exact stability condition in \cite{Lee16} is generally incorrect. We first introduce the classic $L$-index, which is similar in spirit to the condition in \cite{Lee16}. We provide theoretical interpretation and justification of $L$-index. The condition in \cite{Lee16} is then introduced, the similarity between the two indices are drawn and the flaw in \cite{Lee16} is pointed out. A two-bus system counterexample is given to demonstrate the ineffectiveness of the proposed condition in \cite{Lee16}. 





\end{document}